\input amstex
\magnification=\magstep1 
\baselineskip=13pt
\documentstyle{amsppt}
\vsize=8.7truein \CenteredTagsOnSplits \NoRunningHeads
\def\conv{\operatorname{conv}}
\def\PP{\Cal{P}}
\def\BB{\Cal{B}}
\def\CC{\Cal{C}}

\topmatter

\title Neighborliness of the symmetric moment curve 
\endtitle
\author Alexander Barvinok, Seung Jin Lee, and Isabella Novik \endauthor
\address Department of Mathematics, University of Michigan, Ann Arbor,
MI 48109-1043  \endaddress
\email barvinok$\@$umich.edu \endemail
\address Department of Mathematics, University of Michigan, Ann Arbor,
MI 48109-1043 \endaddress
\email lsjin$\@$umich.edu \endemail
\address Department of Mathematics, University of Washington, Seattle, WA 98195-4350 \endaddress
\email novik$\@$math.washington.edu \endemail 
\date November 2011 \enddate
\thanks  The research of the first and second authors was partially supported by NSF Grant DMS-0856640; 
the research of the third author was partially supported by NSF Grants DMS-0801152 and DMS-1069298.
\endthanks 
\abstract We consider the convex hull $\BB_k$ of the symmetric moment curve 
$U_k(t)=\bigl(\cos t, \sin t, \cos 3t, \sin 3t, \ldots, 
\cos (2k-1)t, \sin (2k-1)t \bigr)$ in ${\Bbb R}^{2k}$, 
where $t$ ranges over the unit circle ${\Bbb S}={\Bbb R}/2\pi {\Bbb Z}$. 
The curve $U_k(t)$ is locally neighborly:
as long as $t_1, \ldots, t_k$ lie in an open arc of ${\Bbb S}$ of a certain length $\phi_k>0$, the
convex hull of the points $U_k\left(t_1\right), \ldots, U_k\left(t_k\right)$ is a face of $\BB_k$. 
We characterize the maximum possible length $\phi_k$, proving, in particular, 
that $\phi_k > \pi/2$ for all $k$ and that the limit of $\phi_k$ is $\pi/2$ as $k$ grows. 
This allows us to construct centrally symmetric polytopes with a record number of faces.
\endabstract
\endtopmatter

\document

\head 1. Introduction and main results \endhead

The main object of this paper is the {\it symmetric moment curve} that for a fixed
$k$ lies  in ${\Bbb R}^{2k}$ and is defined by
$$U(t)=U_k(t)=\Bigl(\cos t, \ \sin t, \ \cos 3t, \ \sin 3t, 
\ldots, \cos (2k-1) t, \ \sin (2k-1)t \Bigr).$$
We note that 
$$U(t + \pi) = - U(t) \quad \text{for all} \quad t \in {\Bbb R}.$$
Since $U$ is periodic, we consider $U$ to be defined on the unit 
circle ${\Bbb S} = {\Bbb R}/2 \pi {\Bbb Z}$. 
In particular, for every $t \in {\Bbb S}$, the points $t$
and $t + \pi$ are antipodal points on the circle.

We define the convex body $\BB \subset {\Bbb R}^{2k}$ 
as the convex hull of the symmetric moment curve
$$\BB=\BB_k=\conv\Bigl(U(t): \quad t \in {\Bbb S} \Bigr).$$
Hence $\BB$ is symmetric about the origin, $\BB =-\BB$. 
We note that $\BB_k$ has a non-empty interior 
in ${\Bbb R}^{2k}$ since $U_k(t)$ does not lie in an affine hyperplane.

We are interested in the facial structure of $\BB$ 
(a {\it face} of a convex body is the intersection of the body 
with a supporting affine hyperplane, see, for example, 
Chapter~II of \cite{Ba02}). The convex body $\BB_k$ was 
introduced in \cite{BN08} in the hope that an appropriate 
discretization of $\BB_k$ produces centrally symmetric polytopes with many faces
(although $\BB_2$ was considered first by Smilansky \cite{Sm85} within a certain family
of 4-dimensional convex bodies). Besides being of intrinsic interest, such polytopes appear  
in problems of sparse signal reconstruction, see \cite{Do04},
\cite{RV05}. 
An analogy with the (ordinary) trigonometric moment curve $M(t) \subset {\Bbb R}^{2k}$,
$$M(t)=\Bigl(\cos t, \ \sin t, \ \cos 2t, \ \sin 2t, \ldots, \cos kt,\ \sin kt \Bigr),$$
provided the initial motivation. 
The convex hull of a set of arbitrary $n$ points on the curve $M(t)$ is 
the {\it cyclic polytope}, which was first investigated by Carath\'eodory \cite{Ca11} and 
later by Motzkin \cite{Mo57} and 
Gale \cite{Ga63}.  In the case of an odd dimension $2k+1$ one can define the cyclic polytope as the 
convex hull of $n$ points on the curve $\left(t, t^2, \ldots, t^{2k+1}\right)$.
Cyclic polytopes maximize the number of faces of all dimensions in the class of
polytopes of a given dimension and with a given number of vertices.
This is the famous Upper Bound Theorem conjectured by Motzkin \cite{Mo57} 
and proved by McMullen \cite{Mc70}.
Such maximizers in the class of centrally symmetric polytopes are not known at present.

There are certain similarities between the curves $M(t)$ and $U(t)$. 
An affine hyperplane in ${\Bbb R}^{2k}$ intersects $M(t)$ in no more 
than $2k$ points and it is clear that the bound $2k$ cannot be made smaller. An affine hyperplane in 
${\Bbb R}^{2k}$ intersects $U(t)$ in no more than $4k-2$ points (see Theorem 3.1 below)
 and it is again clear that the bound $4k-2$ 
cannot be made smaller in the class of centrally symmetric curves, 
since a hyperplane passing through the origin 
and some $2k-1$ points on the curve will necessarily intersect the 
curve also in the antipodal $2k-1$ points. 

One crucial feature of the convex hull 
$$\CC=\CC_{k} =\conv\Bigl( M(t): \quad t\in {\Bbb S}\Bigr)$$
of the standard trigonometric moment curve is that it is {\it neighborly}, namely that for any 
$n \leq k$ distinct points $t_1, \ldots, t_n \in {\Bbb S}$, the convex hull 
$\conv\Bigl(M\left(t_1\right), \ldots, M\left(t_n\right)\Bigr)$ is a face of $\CC_k$, see for example,
Chapter II of \cite{Ba02}. One of our main results is that 
the convex hull $\BB_k$ of the symmetric moment curve is neighborly to a large extent.

\proclaim{(1.1) Theorem} For every positive integer $k$ there exists a number 
$${\pi \over 2} \ < \ \phi_k \ < \ \pi$$
such that for an arbitrary open arc $\Gamma \subset {\Bbb S}$ of length $\phi_k$ 
and arbitrary distinct $n \leq k$ points $t_1, \ldots, t_n \in \Gamma$, the set 
$$\conv\Bigl(U\left(t_1\right), \ldots, U\left(t_n\right)\Bigr)$$
is a face of $\BB_k$.
\endproclaim

 It is worth mentioning that Lemma 3.4 below implies that many 
(but not all) $k$-vertex faces of $\BB_k$ are simplices. 
More precisely, if $F$ is a face of $\BB_k$ whose vertex set is $\{U(t_1),\ldots, U(t_k)\}$ 
where $t_1,\ldots,t_k\in {\Bbb S}$ lie in an open semicircle, 
then $F$ is a $(k-1)$-dimensional simplex. Also, since the intersection of faces is a face, it is enough
to verify Theorem 1.1 for $n=k$.

In what follows, we denote by $\phi_k$ the largest possible value that satisfies Theorem~1.1.
We provide a characterization of $\phi_k$, which, in principle, 
allows one to compute it, at least numerically and at least for moderate values of $k$.
This characterization is, roughly, as 
follows. It is not hard to argue that if $\Gamma$ is an open arc of length $\phi_k$ 
then there must be a way to move some of the points $t_1, \ldots, t_k \in \Gamma$ 
towards the endpoints of $\Gamma$, so that the limit position of the affine hyperplane 
supporting $\BB_k$ at $U\left(t_1\right), \ldots, U\left(t_k\right)$ 
will touch the curve $U(t)$ somewhere else. We prove that this limiting configuration 
is as degenerate as it can possibly be: each point $t_i$ collides with one of the endpoints 
$a$ or $b$ of the arc $\Gamma$ so that if $q$ such points collide with $a$
and $k-q$ points with $b$, then the necessarily unique affine hyperplane tangent to 
the symmetric moment curve $U(t)$ at $t=a$ and $t=b$ with multiplicities $2q$ and $2k-2q$, 
respectively, is also tangent to $U(t)$ at some other point (cf.~Theorem 5.1 below).

 We have
$$\phi_2 = {2 \pi \over 3} \approx 2.094395103,$$
(this follows from results of Smilansky \cite{Sm85}), and we computed 
$$\phi_3 = \pi - \arccos {3 - \sqrt{5} \over 2} \approx 1.962719003$$
as well as 
$$\split \phi_4 =&2 \arccos\left( -{1 \over 48} \left(91 + 336 \sqrt{15}\right)^{1/3} + 
{119 \over 48\left(91 +336 \sqrt{15}\right)^{1/3}}
+{29 \over 48}\right)\\
\approx &1.870658532,\endsplit$$
cf.~Example 5.2.

It is worth noting that in \cite{BN08} we were only able to 
verify that $\phi_k >0$, while in \cite{Le11} 
the first explicit lower bound $\phi_k > \sqrt{6} k^{-3/2}$ was established.

We conjecture that for an even $k$ the value of $\phi_k$ in Theorem 1.1 
satisfies $\phi_k =2 \alpha_k$, where $\alpha_k$ is the smallest positive root of the equation 
$$\cos \alpha +1 +\sum_{j=1}^{k-1} (-1)^j {(2j-1)!! \over (2j)!!} \tan^{2j} \alpha =0.$$
Here $n!!$ is the product of positive integers not exceeding $n$ and having the same 
parity as $n$. Some supporting evidence for the conjecture is provided in Section 7.

Theorem 1.1 allows one to construct $2k$-dimensional centrally symmetric polytopes with $n$ vertices and 
$$ \left(4 \cdot 2^{-k} -4^{-k+1} +O\left({1 \over n}\right)\right){n \choose k}$$ faces of 
dimension $k-1$, a new record. Indeed, suppose that $n=4m$ and consider a centrally 
symmetric configuration $A=A_0 \cup A_1 \cup A_2 \cup A_3$
of $4m$ distinct points in ${\Bbb S}$, where each set $A_j$ 
contains $m$ distinct points in the vicinity of $j\pi/2$ for $j=0, 1, 2, 3$ 
so that the length of any arc whose endpoints are 
in $A_j$ and $A_{(j+1) \mod 4}$ is less than $\phi_k$. 
Letting 
$$P=\conv\Bigl( U(t): \quad t \in A\Bigr),$$
we observe that $P$ is a centrally symmetric polytope with $4m$ vertices and at least 
$4 {2m \choose k}-4{m \choose k}$ faces of dimension $k-1$.  
For example, if $k=2$ we obtain a 4-dimensional 
centrally symmetric polytope with $4m$ vertices such that 
approximately $3/4$ of all pairs of vertices are guaranteed to span edges of the 
polytope. Curiously, if we consider instead the convex hull of 
$4m$ points $U_2\left(t_i\right)$ for points $t_i$ 
uniformly distributed on the circle ${\Bbb S}$ then only about 
$2/3$ of all pairs of vertices span edges of the resulting 
4-dimensional centrally symmetric polytope, see \cite{BN08}.

In \cite{B+11}, we apply Theorem 1.1 to construct various 
families of centrally symmetric polytopes with many faces.
Namely, we construct a $d$-dimensional centrally symmetric polytope $P$ with about 
$3^{d/4} \approx (1.316)^d$ vertices such that 
every pair of non-antipodal vertices of $P$ spans an edge of $P$. 
For $k \geq 1$, we construct a $d$-dimensional centrally symmetric polytope $P$ 
with an arbitrarily large number $n$ of vertices such that 
the number of $k$-dimensional faces of $P$ is at least 
$(1-(\delta_k)^d){n \choose k+1}$ for some
$0 < \delta_k < 1$ (we have $\delta_1 \approx 3^{-1/4} 
\approx 0.77$ and $\delta_k \approx \left(1-5^{-k}\right)^{1/(6k+4)}$
for $k >1$). Finally, for an integer $k > 1$ and $\alpha>0$, 
we construct a centrally symmetric polytope $P$ with an arbitrarily large 
number $n$ of vertices and of dimension $d=k^{1+o(1)}$ such 
that the number of $k$-dimensional faces of 
$P$ is at least $\left(1-k^{-\alpha}\right) {n \choose k+1}$.

 It is important to note that we are interested in asymptotics of 
the number of faces when the dimension is fixed and the number 
of vertices grows, in contrast to a number of 
results in the literature where both the dimension and the number 
of vertices grow, see \cite{RV05}, \cite{LN06},
and \cite{DT09}.

We prove that $\pi/2$ is the limit of neighborliness of $U_k(t)$ as $k$ grows.

\proclaim{(1.2) Theorem} Let $\phi_k$ be the largest number satisfying Theorem 1.1. Then 
$$\lim_{k \longrightarrow +\infty} \phi_k = {\pi \over 2} \approx 1.570796327.$$
\endproclaim

Our computations strongly suggest that the values $\phi_k$ are monotone decreasing, but 
we were unable to prove that.  

The complete facial structure of $\BB_2$ was described by Smilansky \cite{Sm85}: 
the 0-dimensional faces are the points $U(t)$ for $t \in {\Bbb S}$, 
the 1-dimensional faces are the intervals $[U(a), U(b)]$, where $a, b \in {\Bbb S}$ and
the length of the shorter arc with the endpoints $a$ and $b$ is less than $2\pi/3$, 
and the 2-dimensional faces are the triangles with the vertices $\bigl\{U(a), U(b), U(c)\bigr\}$, 
where $a, b, c \in {\Bbb S}$ are vertices of an equilateral triangle.
There are no other faces of $\BB_2$. 

In \cite{BN08} and \cite{Vi11} the edges (1-dimensional faces) 
of $\BB_{k}$ were completely characterized. 
Namely,  it was shown in \cite{BN08} that for $a\ne b$ the interval $\left[ U(a),\ U(b) \right]$ 
is an edge of $\BB_{k}$ if the length of the 
shorter arc with the endpoints $a,  b \in {\Bbb S}$ 
is smaller than ${2k-2 \over 2k-1} \pi$, and 
it was shown in \cite{Vi11} that there are no other edges. 

Already for $\BB_3$, the complete facial structure is not known; 
some faces of $\BB_3$ were computed in \cite{Li07}.

We obtain the following partial result showing some``connectedness'' of faces of $\BB_k$.
\proclaim{(1.3) Theorem} Let $\Gamma \subset {\Bbb S}$ be an open arc with the 
endpoints $a$ and $b$ and let $\overline{\Gamma}$ be the closure of $\Gamma$. Let  
$t_2, \ldots, t_k \in {\Bbb S} \setminus  \overline{\Gamma}$
be distinct points such that the set $\overline{\Gamma} \cup \bigl\{t_2, \ldots, t_k \bigr\}$ 
lies in an open semicircle in ${\Bbb S}$.
Suppose that the points $U(a), U(t_2), \ldots, U(t_k)$ lie in a face of $\BB_k$ and that the points 
$U(b), U(t_2), \ldots, U(t_k)$ lie in a face of $\BB_k$. Then for all $t_1 \in \Gamma$  the set
$$\conv\Bigl(U(t_1), \ldots, U(t_k)\Bigr)$$
is a face of $\BB_k$.
\endproclaim

\example{(1.4) Example} Let $k=3$ and let $t_1, t_2, t_3 \in {\Bbb S}$ 
be any points such that the length of the arc with the endpoints $t_1$ and $t_2$ 
is less than $2\pi/5$, the length of the arc with the endpoints $t_2$ and $t_3$
is less than $2\pi/5$, and $t_2$ lies between $t_1$ and $t_3$. Then 
$$\conv\Bigl( U\left(t_1\right), U\left(t_2\right), U\left(t_3\right) \Bigr) \tag1.4.1$$
is a face of $\BB_3$. Indeed, without loss of generality, we may assume that $t_2=0$.  
If  $t_1=2\pi/5$ and $t_3=-2\pi/5$, the triangle (1.4.1) lies in the  face of 
$\BB_3$ determined by the equation $\cos 5t =1$. If we move $t_3$ sufficiently 
close to $t_2=0$, then (1.4.1) is a face from our estimate of $\phi_3$ in Theorem 1.1. Therefore, by 
Theorem 1.3, for $t_1=2\pi/5$, $t_2=0$, 
and all $-2\pi/5 < t_3 < 0$, the set (1.4.1) is a face of $\BB_3$.
Let us now fix some $-2\pi/5 < t_3 < 0$. If we move $t_1$ sufficiently close to $t_2=0$ then 
(1.4.1) is a face by Theorem 1.1. Applying Theorem 1.3 again we conclude that for all 
$0 < t_1 < 2\pi/5$, $-2\pi/5 < t_3 < 0$, and $t_2=0$, the set (1.4.1) is a face of $\BB_3$.
\endexample

In the rest of the paper, we prove Theorems 1.1--1.3. 
In Section 2, we introduce the analytic language of 
raked trigonometric polynomials which supplants the 
geometric language of affine hyperplanes in our proofs. 
We also outline the plan for the proofs.

\head 2. Polynomials.  The plan of the proofs \endhead 

\subhead (2.1) Raked trigonometric and complex polynomials \endsubhead 
We consider {\it raked trigonometric polynomials} of degree at most $2k-1$:
$$f(t) =c + \sum_{j=1}^k a_j \cos (2j-1) t + \sum_{j=1}^k b_j \sin (2j-1) t \quad \text{for} \quad 
t \in {\Bbb S}, \tag2.1.1$$
where $c, a_j, b_j \in {\Bbb R}$. We say that $\deg f =2k-1$ if $a_k \ne 0$ or if $b_k \ne 0$.
 Equivalently, we can write 
$$f(t) = c + \bigl\langle C,\ U(t) \big\rangle,$$
where $C=\left(a_1, b_1, \ldots, a_k, b_k\right) \in {\Bbb R}^{2k}$ and 
$\langle \cdot, \cdot \rangle$ is the standard scalar product in ${\Bbb R}^{2k}$.

Writing 
$$\cos nt = {e^{int} + e^{-int} \over 2} \quad \text{and}\quad  \sin nt ={e^{int} - e^{-int} \over 2i}$$
and substituting $z=e^{it}$, we associate with (2.1.1) a complex polynomial 
$$\PP(f)(z)=z^{2k-1} \left( c + \sum_{j=1}^k a_j {z^{2j-1} + z^{1-2j} \over 2} 
+\sum_{j=1}^k b_j {z^{2j-1} -z^{1-2j} \over 2i} \right). \tag2.1.2$$
Hence 
$$\deg \PP(f) \ \leq \ 4k-2. \tag2.1.3$$ 
Moreover, if $\deg f =2k-1$ then for $p=\PP(f)$ we have $\deg p =4k-2$ and $p(0) \ne 0$.
Since 
$$\cos(t+a)=\cos t \cos a - \sin t \sin a \quad \text{and} 
\quad \sin(t+a)=\sin t \cos a +\cos t \sin a,$$
for any fixed $a \in {\Bbb S}$ and any raked trigonometric polynomial $f(t)$, the function 
$$h(t)=f(t+a) \quad \text{for} \quad t \in {\Bbb S}$$ 
is also a raked trigonometric polynomial of the same degree.

\definition{(2.2) Definition} We say that a point $t^{\ast} \in {\Bbb S}$ is a 
{\it root of multiplicity} $m$ (where $m \geq 1$ is 
an integer) of a trigonometric polynomial $f$, if 
$$f\left(t^{\ast}\right) = \ldots = f^{(m-1)}\left(t^{\ast}\right) =0$$
and 
$$f^{(m)}\left(t^{\ast}\right) \ne 0.$$
Similarly, we say that a number $z^{\ast} \in {\Bbb C}$ is a root of multiplicity $m$ 
of a polynomial $p(z)$ if 
$$p\left(z^{\ast}\right)= \ldots =p^{(m-1)}\left(z^{\ast}\right)=0$$ and 
$$p^{(m)}\left(z^{\ast}\right)\ne 0.$$
\enddefinition

\remark{(2.3) Remark}
\roster
\item We note that 
$$\conv\Bigl(U\left(t_1\right), \ldots, U\left(t_n\right) \Bigr)$$
for distinct $t_1, \ldots, t_n \in {\Bbb S}$
is a face of $\BB_k$ if and only if there exists a raked 
trigonometric polynomial $f(t)\not\equiv 0$ of
degree at most $2k-1$ such that (i) each $t_j$, 
$j=1,\ldots, n$, is a root of $f$ of an even multiplicity,
and (ii) $f$ has no other roots. The roots of $f$ should have 
even multiplicities as for $f$ to determine 
a face, the values of $f$ must not change sign on ${\Bbb S}$. 

\item We will often use the following observation: 
if $f$ is a trigonometric polynomial with constant term 1 
that does not change sign on ${\Bbb S}$ then $f(t) \geq 0$ for all $t \in {\Bbb S}$, since 
$${1 \over 2\pi} \int_{\Bbb S} f(t) \ dt =1.$$
\endroster
\endremark

\subhead (2.4) Plan of the proofs \endsubhead 

In Section 3, we prove basic facts about the roots of 
raked trigonometric polynomials. We bound their number 
and restrict possible positions in the circle ${\Bbb S}$ 
(Theorem 3.1). While the total number of roots of 
a raked trigonometric polynomial $f$ with $\deg f \leq 2k-1$ 
does not exceed $4k-2$, we prove that 
the number of roots, counting multiplicities, in any open arc 
$\Gamma \subset {\Bbb S}$ of length less than $\pi$
may not exceed $2k$; moreover, if that number is equal to $2k$
then any additional root of $f$ has to lie in 
the opposite arc $\Gamma +\pi$.

In Section 4, we consider one-parameter families $f_s(t)$ 
of raked trigonometric polynomials with constant term 1,
obtained by fixing some roots of $f_s(t)$ and moving one, possibly multiple, root,
so that the total number of controlled roots is $2k$, counting multiplicities, and the roots 
remain in an open semi-circle of ${\Bbb S}$. We prove that 
$(\partial/\partial s) f_s(t) \not\equiv 0$ (Theorem 4.2).
In geometric terms, Theorem 4.2 implies that if we choose $k$ 
distinct points $t_1, \ldots, t_k$ lying in 
an open semicircle of ${\Bbb S}$, consider the (necessarily unique) affine hyperplane $H$ tangent 
to the symmetric moment curve $U_k(t)$ at $t_1, \ldots, t_k$ and then start moving the point
$t_1$, while keeping the points $t_2, \ldots, t_n$ intact, 
the velocity of the unit normal of $H$ is never zero. 

In Section 5, we prove a characterization (see Theorem 5.1 and Lemma 5.7) 
of the value of $\phi_k$ introduced 
in Theorem 1.1. In analytic terms, we prove that if $\Gamma \subset {\Bbb S}$ 
is an open arc of a certain length and if a raked trigonometric polynomial $f(t)$ 
has $2k$ roots, counting multiplicities, 
in $\Gamma$, then $f$ has no other roots in ${\Bbb S}$. 
Moreover, we prove that for the maximum possible 
length $\phi_k$ of such an arc $\Gamma$, there are positive 
even integers $m_a$ and $m_b$ such that 
$m_a+m_b=2k$ and such that the unique, up to a non-zero multiple, 
raked trigonometric polynomial $f(t)$ of 
degree $2k-1$ that has a root of multiplicity $m_a$ at one endpoint  
of $\Gamma$ and a root of multiplicity $m_b$ at the other endpoint  
of $\Gamma$ also has a root of an even multiplicity in 
$\Gamma +\pi$ and does not change its sign on ${\Bbb S}$.

In Section 6, we prove Theorems 1.1 and 1.3. In particular, we prove that for every 
positive integer $k$ there exists a number $\phi_k > \pi/2$ with the following property:
 if $f(t)$ is an arbitrary raked trigonometric polynomial of degree $2k-1$,  with constant term 1,
and such that $f(t)$ has $2k$ roots, counting multiplicities, in an open arc $\Gamma \subset {\Bbb S}$ 
of length $\phi_k$ and all roots in $\Gamma$ have even multiplicities, 
then $f(t)$ is positive everywhere else in ${\Bbb S}$ (Theorem~6.1).

In Section 7, we prove Theorem 1.2. Also, for an even $k$, 
we deduce an equation for the value of $0< \alpha < \pi/2$ 
such that the unique raked trigonometric polynomial of degree 
$2k-1$ with constant term 1 that has roots at $t=\pm \alpha$ 
of multiplicity $k$ each also has a root of an even multiplicity 
at $t=\pi$ while remaining non-negative on ${\Bbb S}$.
We conjecture that $\phi_k=2 \alpha$.

\head 3. Roots and multiplicities \endhead 

We consider raked trigonometric polynomials $f(t)$ defined by (2.1.1).
In this section we prove the following main result.

\proclaim{(3.1) Theorem} 
Let $f(t)\not\equiv 0$ be a raked trigonometric polynomial of degree at most $2k-1$,
let $t_1, \ldots, t_n \in {\Bbb S}$ be distinct roots of $f$ in ${\Bbb S}$, and let 
$m_1, \ldots, m_n$ be their multiplicities.
\roster
\item We have
$$\sum_{i=1}^n m_i \ \leq \ 4k-2.$$
\item If the constant term of $f$ is 0 and the set 
$\left\{t_1, \ldots, t_n\right\}$ does not contain a 
pair of antipodal points, then
$$\sum_{i=1}^n m_i \ \leq \ 2k-1.$$
\item If $t_1, \ldots, t_n$ lie in an open semicircle of ${\Bbb S}$, then 
$$\sum_{i=1}^n m_i \ \leq \ 2k.$$
\item Suppose that $t_1, \ldots, t_n$ lie in an arc 
$\Gamma \subset {\Bbb S}$ of length less than $\pi$,
  that
$$\sum_{i=1}^n m_i =2k,$$
and that $t^{\ast} \in {\Bbb S} \setminus \Gamma$ is yet another root of $f$. Then 
$t^{\ast} \in \Gamma +\pi$.
\endroster
\endproclaim

To prove Theorem 3.1, we establish a correspondence between 
the roots of a trigonometric polynomial $f(t)$ 
and those of the corresponding complex polynomial $p(z)=\PP(f)$ defined by (2.1.2).

\proclaim{(3.2) Lemma} 
A point $t^{\ast} \in {\Bbb S}$ is a root of multiplicity $m$ of $f(t)$ if and only if
 $z^{\ast} =e^{i t^{\ast}}$ is a root of multiplicity $m$ of $\PP(f)$.
\endproclaim
\demo{Proof} Let $p=\PP(f)$.
It follows from (2.1.2) that
$$p\left(e^{it}\right)=e^{(2k-1)it}f(t). \tag3.2.1$$
Differentiating (3.2.1), we infer by induction
 that
$$ i^r\sum_{j=1}^r d_{j,r} e^{ijt} p^{(j)}\left(e^{it}\right)
=\sum_{j=0}^r i^{r-j}  c_{j,r} e^{(2k-1)it}\cdot f^{(j)}(t)
\quad \text{for all} \quad r\geq 1,$$
where the constants $c_{j,r}$, $d_{j,r}$ are positive integers. Thus, if
$f^{(r)}(t^{\ast})$ is zero for $r=0,1,\ldots,m-1$ and nonzero for $r=m$,
then so is $p^{(r)}\left(e^{it^{\ast}}\right)$, and vice versa. The statement now follows. 
{\hfill \hfill \hfill}\qed
\enddemo

\subhead (3.3) Proof of Theorem 3.1 \endsubhead 

Part (1) follows from Lemma 3.2 and bound (2.1.3). 

If $f$ has a zero constant term, then $f$ satisfies 
$$f(t+\pi) =-f(t) \quad \text{for all} \quad t \in {\Bbb S}.$$
Then  $t_i +\pi$ is a root of $f(t)$ of multiplicity $m_i$ and the proof 
of Part (2) follows from Part (1).

To prove Part (3), let $g(t)=f'(t)$. Then $g$ has a zero constant term.
If 
$$\sum_{i=1}^n m_i > 2k,$$
then by Rolle's Theorem, the total number of roots of $g(t)$ in the semicircle, 
counting multiplicities, is at least $2k$, and so $g(t) \equiv 0$ by Part (2), 
which is a contradiction.

To prove Part (4), we assume without loss of generality that $t_1, \ldots, t_n$ is the order of 
the roots on the arc $\Gamma$ and let $\tilde{\Gamma}$ 
be the closed arc with the endpoints $t_1$ and $t_n$. 
By Rolle's Theorem, the total number of roots of $g(t)$, counting multiplicities, in $\tilde{\Gamma}$ 
is at least $2k-1$, and hence the total number of roots of $g(t)$, counting multiplicities, in 
$\tilde{\Gamma} \cup\left(\tilde{\Gamma}+\pi\right)$ is at least $4k-2$.
 If $t^{\ast} \notin \tilde{\Gamma} \cup \left(\tilde{\Gamma} +\pi\right)$, then by Rolle's theorem 
there is a root of $g(t)$ outside of $\tilde{\Gamma} \cup \left(\tilde{\Gamma} +\pi\right)$,
 and hence the total number of roots of $g(t)$ in  ${\Bbb S}$, counting multiplicities, 
is at least $4k-1$. Thus, by Part (1), $g(t) \equiv 0$, which is a contradiction.
{\hfill \hfill \hfill} \qed

We will utilize the following geometric corollary of Theorem 3.1.

\proclaim{(3.4) Lemma} 
Let $t_1, \ldots, t_n \in {\Bbb S}$ be distinct points lying in an open semicircle and let
$m_1, \ldots, m_n$ be positive integers such that 
$$\sum_{i=1}^n m_i =2k.$$
Then the following $2k$ vectors
$$\split &U\left(t_i\right)-U\left(t_n\right) \quad \text{for} \quad i=1, \ldots, n-1, \\
&{d^j \over dt^j} U(t) \Big|_{t=t_i} \quad \text{for} 
\quad j=1, \ldots,  m_i-1 \quad \text{if} \quad m_i>1 \quad \text{and} \quad i=1, \ldots, n, \\
&{d^{m_1} \over dt^{m_1}} U(t) \Big|_{t=t_1} \endsplit$$
are linearly independent in ${\Bbb R}^{2k}$. 
\endproclaim
\demo{Proof} 
Seeking a contradiction, we assume that the vectors are not linearly independent. Then there exists 
a non-zero vector $C \in {\Bbb R}^{2k}$ orthogonal to all these $2k$ vectors. 
Consider the raked trigonometric polynomial 
$$f(t)=\langle C,\ U(t) -U\left(t_n\right) \rangle \quad \text{for} \quad t \in {\Bbb S}.$$
Then $t_1, \ldots, t_n$ are roots of $f(t)$. Moreover, 
the multiplicity of $t_i$ is at least $m_i$ for $i>1$ and at least $m_1+1$ 
for $i=1$. It follows from Part (3) of Theorem 3.1 
that $f(t) \equiv 0$, which contradicts that $C \ne 0$.
{\hfill \hfill \hfill} \qed
\enddemo

Finally, we prove that a raked trigonometric polynomial is determined, up to a constant factor,
 by its roots of the total multiplicity $2k$ provided those roots lie in an open semicircle.

\proclaim{(3.5) Corollary} 
Let $t_1, \ldots, t_n \in {\Bbb S}$ be distinct points lying in an open semicircle, and let
$m_1, \ldots, m_n$ be positive integers such that 
$$\sum_{i=1}^n m_i =2k.$$
Then there exists a unique raked trigonometric polynomial $f(t)$ of degree at most $2k-1$ 
and with constant term 1, such that $t_i$ is a root of $f(t)$ of multiplicity $m_i$ for 
all $i=1, \ldots, n$. 
Moreover, $f$ depends analytically on $t_1, \ldots, t_n$.
\endproclaim 
\demo{Proof} Such a polynomial $f(t)$ can be written 
as
$$f(t)= \langle C,\ U(t)-U(t_n) \rangle  \quad \text{for} \quad t \in {\Bbb S},$$
where $C  \in {\Bbb R}^{2k}$ is orthogonal to the $2k-1$ vectors 
$$\split &U\left(t_i\right)-U\left(t_n\right) \quad \text{for} \quad i=1, \ldots, n-1, \\
&{d^j \over dt^j} U(t) \Big|_{t=t_i} \quad \text{for} \quad j=1, \ldots,  
m_i-1 \quad \text{if} \quad m_i>1 \quad \text{and} \quad i=1, \ldots, n. \endsplit$$
By Lemma 3.4, these $2k-1$ vectors span a hyperplane in ${\Bbb R}^{2k-1}$ and hence,
up to a scalar, there is a unique choice of $C$. By Part (2) of Theorem 3.1, 
$f$ has a non-zero constant term if $C \ne 0$. Therefore, there is a 
unique choice of $C$ that makes the constant term of $f(t)$ equal 1.
By Part (3) of Theorem 3.1 the multiplicities of the roots $t_i$ are exactly $m_i$ for 
$i=1, \ldots, n$. 
{\hfill \hfill \hfill} \qed
\enddemo
Note that in fact $\deg f=2k-1$. This follows from Part (3) of Theorem 3.1.

We will also need the following ``deformation construction''.

\proclaim{(3.6) Lemma} 
Let $f(t)$ be a raked trigonometric polynomial of degree $2k-1$ such that $f(-t)=f(t)$
for all $t\in{\Bbb S}$,
and let $p=\PP(f)$ be the corresponding complex polynomial associated with $f$ via (2.1.2).
Then $p(0)\ne 0$ and the multiset $M$ of roots of $p$ can be split into $2k-1$ unordered pairs
$\left\{\zeta_j, \zeta_j^{-1}\right\}$ for $j=1, \ldots, 2k-1$. 
Moreover, for any real $\lambda \ne 0$, the multiset 
$M_{\lambda}$ consisting of $2k-1$ unordered pairs $\left\{ \xi_j, \xi_j^{-1}\right\}$
defined by
$$\xi_j + \xi_j^{-1} = \lambda \left(\zeta_j + \zeta_j^{-1}\right) 
\quad \text{for} \quad j=1, \ldots, 2k-1$$
is the multiset of roots of a certain complex polynomial $p_{\lambda}$ such that 
$p_{\lambda} =\PP\left(f_{\lambda}\right)$ for a raked trigonometric polynomial 
$f_{\lambda}(t)$ of degree $2k-1$ satisfying $f_{\lambda}(-t)=f_{\lambda}(t)$.
\endproclaim
\demo{Proof}
This is Lemma 5.1 of \cite{BN08}.
{\hfill \hfill \hfill} \qed
\enddemo

We call $f_{\lambda}(t)$ a $\lambda$-{\it deformation} of $f$.

\head 4. Parametric families of trigonometric polynomials \endhead 

\subhead (4.1) Parametric polynomials \endsubhead 
Let $\Gamma \subset {\Bbb S}$ be an open arc. We consider 
raked trigonometric polynomials 
$$f_s(t) = 1+ \sum_{j=1}^k a_j(s) \cos (2j-1) t 
+\sum_{j=1}^k b_j(s) \sin (2j-1) t \quad \text{for} \quad 
t \in {\Bbb S}, \tag4.1.1$$
where $a_j(s)$ and $b_j(s)$ are real analytic functions of $s \in \Gamma$.
We define 
$$g_s(t)={\partial \over \partial s} f_s(t),$$
and so
$$g_s(t)=\sum_{j=1}^k a_j'(s) \cos (2j-1) t + \sum_{j=1}^k b_j'(s) \sin (2j-1) t. \tag4.1.2$$

The goal of  this section is to prove the following result.
\proclaim{(4.2) Theorem} 
Let $\Gamma \subset {\Bbb S}$ be an open arc, let 
$t_2, \ldots, t_n \in {\Bbb S} \setminus \Gamma$ be distinct 
points such that the set $\Gamma \cup \left\{t_2, \ldots, t_n \right\}$ lies in an open semicircle, and let 
$m_1, \ldots, m_n $ be positive integers such that 
$$\sum_{i=1}^n m_i =2k.$$
For every $s \in \Gamma$, let $f_s(t)$ be the unique raked trigonometric polynomial of degree $2k-1$ 
with constant term 1 such that 
for $i=2, \ldots, n$ the point $t_i$ is a root of $f_s(t)$ of multiplicity $m_i$ and 
$s$ is a root of $f_s(t)$ of multiplicity $m_1$, cf.~Corollary 3.5. Define 
$$g_s(t) = {\partial \over \partial s} f_s(t).$$
Then 
$$g_s(t) \not\equiv 0 \quad \text{for all} \quad s \in \Gamma.$$
\endproclaim

To prove Theorem 4.2, we use the notion of the wedge product.

\subhead (4.3) Wedge product \endsubhead 
Given linearly independent vectors $V_1, \ldots, V_{2k-1} \in {\Bbb R}^{2k}$
we define their wedge product 
$$W=V_1 \wedge \ldots \wedge V_{2k-1}$$ 
as the unique vector $W$ orthogonal to the hyperplane spanned 
by $V_1, \ldots, V_{2k-1}$ whose length is the 
volume of the $(2k-1)$-dimensional parallelepiped spanned by 
$V_1, \ldots, V_{2k-1}$ and such that the basis 
$V_1, \ldots, V_{2k-1}, W$ is co-oriented with the standard basis of ${\Bbb R}^{2k}$. 
If vectors $V_1, \ldots, V_{2k-1}$ are linearly dependent, we let
$$V_1 \wedge \ldots \wedge V_{2k-1} =0.$$
Suppose that vectors $V_1(s), \ldots, V_{2k-1}(s)$ depend smoothly 
on a real parameter $s$. We will use the following standard fact:
$$\aligned &{d \over ds} \Bigl( V_1(s) \wedge \ldots \wedge V_{2k-1}(s) \Bigr) \\&\quad =
\sum_{j=1}^{2k-1} V_1(s) 
\wedge \ldots \wedge V_{j-1}(s) \wedge {d \over ds} V_j(s) \wedge V_{j+1}(s) \wedge 
\ldots \wedge V_{2k-1}(s). \endaligned \tag4.3.1$$

\subhead (4.4) Proof of Theorem 4.2 \endsubhead 
For $s \in \Gamma$,  consider the following ordered set of $2k-1$ vectors:
$$\aligned &U\left(t_i\right)-U\left(t_n\right) \quad \text{for} \quad i=2, \ldots, n-1, \\
&{d^j \over dt^j} U(t) \Big|_{t=t_i} \quad \text{for} \quad 
j=1, \ldots,  m_i-1 \quad \text{if} \quad m_i>1
\quad  \text{and} \quad i=2, \ldots, n,\\
&U(s)-U\left(t_n\right), \\
&{d^j \over dt^j} U(t) \Big|_{t=s} \quad \text{for} \quad j=1, \ldots, m_1-1 
\quad \text{if} \quad m_1>1.  \endaligned  \tag4.4.1$$
Let $C(s)$ be the wedge product of  vectors of (4.4.1). 
By Lemma 3.4, the vectors of~(4.4.1) are linearly independent 
for all $s \in \Gamma$, and hence $C(s) \ne 0$ for all $s \in \Gamma$. 

For $s \in \Gamma$,  define a raked trigonometric polynomial 
$$F_s(t) =\big\langle C(s), \ U(t) - U(t_n) \big\rangle. \tag4.4.2$$
We note that $F_s(t) \not\equiv 0$ for all $s \in \Gamma$. 
For $i=2, \ldots, n$, the point $t_i$ is a root of $F_s(t)$ 
of multiplicity at least $m_i$ and $s$ is a root of $F_s(t)$ 
of multiplicity at least $m_1$. By Part (3) of Theorem 3.1 the multiplicities are exactly 
$m_i$. Let $\alpha(s)$ be the constant term of $F_s(t)$. Then 
$$\alpha(s) = - \big\langle C(s), U(t_n) \big\rangle.$$
By Part (2) of Theorem 3.1
$$\alpha(s) \ne 0 \quad \text{for all}  \quad s \in \Gamma.$$
Therefore, 
$$f_s(t) = {F_s(t) \over \alpha(s)}.$$
Seeking a contradiction, let us assume that $g_{s}(t) \equiv 0$ for some $s \in \Gamma$.
We have 
$$g_s(t) = {\partial \over \partial s} f_s(t) = 
{\alpha(s) {\partial \over \partial s} F_s(t) - \alpha'(s) F_s(t) \over \alpha^2(s)}.$$ 
If $g_{s}(t) \equiv 0$, then
$$\alpha(s) {\partial \over \partial s} F_s(t) -\alpha'(s) F_{s}(t) \equiv 0,$$
and (4.4.2) yields that 
$$\alpha(s) C'(s) - \alpha'(s) C(s) =0 \tag4.4.3$$
for some $s \in \Gamma$.
Let us consider $C'(s)$, the derivative of the wedge product of (4.4.1). 
Applying formula (4.3.1) we note that all of the $2k-1$ terms of (4.3.1) except the last one are zeros 
since the corresponding wedge product either contains a zero vector or two identical vectors. 
Hence $C'(s)$ is the wedge product of the following ordered set of vectors 
$$\aligned &U\left(t_i\right)-U\left(t_n\right) \quad \text{for} \quad i=2, \ldots, n-1, \\
&{d^j \over dt^j} U(t) \Big|_{t=t_i} \quad \text{for} \quad j=1, 
\ldots,  m_i-1  \quad \text{if} \quad m_i > 1\quad \text{and} \quad i=2, \ldots, n,\\
&U(s)-U\left(t_n\right), \\
&{d^j \over dt^j} U(t) \Big|_{t=s} \quad \text{for}  \quad j=1, \ldots, m_1-2 \quad \text{if} \quad 
m_1>2 \quad \text{and} \quad j=m_1.  \endaligned  \tag4.4.4$$
The wedge products (4.4.1) for $C(s)$ and (4.4.4) for $C'(s)$ differ in two vectors, 
$$A(s)={d^{m_1 -1} \over d t^{m_1-1}} U(t) \Big|_{t=s} \quad \text{and} \quad 
B(s)={d^{m_1} \over d t^{m_1}} U(t) \Big|_{t=s} \quad \text{if} \quad m_1>1$$
and 
$$A(s)=U(s) -U\left(t_n\right) \quad \text{and} \quad B(s)={d \over dt} U(t) \Big|_{t=s} \quad 
\text{if} \quad m_1=1.$$
Vector $A(s)$ is present in (4.4.1) and absent in (4.4.4) while 
vector $B(s)$ is absent in (4.4.1) and present 
in (4.4.4). Therefore, (4.4.3) implies that the set consisting of the vector 
$$\alpha(s) A(s) - \alpha'(s) B(s)$$
and the $2k-2$ vectors common to wedges (4.4.1) and (4.4.4) is linearly dependent.
However, as $\alpha(s) \ne 0$, this contradicts Lemma 3.4.
{\hfill \hfill \hfill} \qed

We will need the following result.

\proclaim{(4.5) Lemma} Let $f_s(t)$ and $g_s(t)$ be trigonometric polynomials (4.1.1) and (4.1.2) 
respectively and let $m$ be a positive integer.
\roster
\item If $t^{\ast} \in {\Bbb S}$ is a root of $f_s(t)$ 
of multiplicity at least $m$ for all $s \in \Gamma$,
then $t^{\ast}$ is a root of $g_s(t)$ of multiplicity at least $m$ for all $s \in \Gamma$.
\item If $m>1$ and $s$ is a root of $f_s(t)$ of multiplicity at least $m$ for all $s \in \Gamma$,
then $s$ is a root of $g_s(t)$ of multiplicity at least $m-1$ for all $s \in \Gamma$.
\endroster
\endproclaim
\demo{Proof} Suppose that
$$f_s\left(t^{\ast}\right) = \ldots = 
{\partial^{m-1} \over \partial t^{m-1}} f_s(t) \Big|_{t=t^{\ast}}=0.$$
Differentiating with respect to $s$ yields Part (1).

Suppose that 
$$f_s(s) = {\partial^j \over \partial t^{j}} f_s(t) \Big|_{t=s} =0 
\quad \text{for} \quad j=1, \ldots, m-1.$$
Differentiating with respect to $s$ we obtain 
$$\split 0 =& {\partial \over \partial s} f_s(t) \Big|_{t=s} 
+ {\partial \over \partial t} f_s(t) \Big|_{t=s} = 
{\partial \over \partial s} {\partial^j \over \partial t^{j}} f_s(t) \Big|_{t=s} 
+ {\partial^{j+1} \over \partial t^{j+1}}
f_s(t) \Big|_{t=s} \\&\qquad \text{for} \quad j=1, \ldots, m-1. \endsplit$$
Therefore,
$$g_s(s) = {\partial^j \over \partial t^j} g_s(t) \Big|_{t=s} =0 \quad 
\text{for} \quad j=1, \ldots, m-2,$$
and the proof of Part (2) follows.
{\hfill \hfill \hfill} \qed
\enddemo 

\head 5. Critical arcs \endhead 

This section is devoted to verifying the following result.
\proclaim{(5.1) Theorem} 
\roster
\item For every $k \geq 1$ there exists a non-empty open arc 
$\Gamma \subset {\Bbb S}$ with the following property: if 
$t_1, \ldots, t_n \in \Gamma$ are distinct points and $m_1, \ldots, m_n$ are positive 
even integers satisfying 
$$\sum_{i=1}^n m_i = 2k,$$
then the unique raked trigonometric polynomial $f(t)$ of degree $2k-1$ 
with constant term 1 that has each point
$t_i$ as a root of multiplicity $m_i$, has no other roots in ${\Bbb S}$. Moreover,  
$f(t) \geq 0$ for all $t \in {\Bbb S}$.
\item Let $\Gamma \subset {\Bbb S}$ be an open arc as in Part (1) of the maximum possible 
length and let $a$ and $b$ be the endpoints of $\Gamma$. Then there are positive 
even integers $m_a$ and $m_b$ such that $m_a + m_b =2k$ and such that 
the unique raked trigonometric polynomial $f(t)$ of degree $2k-1$ with constant term 1 that has a root 
at $t=a$ of multiplicity $m_a$ and a root at $t=b$ of multiplicity $m_b$ 
is non-negative on ${\Bbb S}$ and has 
a root (of necessarily even multiplicity) in the arc $\Gamma + \pi$. 
\item Fix positive even integers $m_a$ and $m_b$ such that $m_a+m_b=2k$. 
Let $\Gamma \subset {\Bbb S}$ be an open 
arc of length less than $\pi$ and let $a$ be an endpoint of $\Gamma$. 
For $b \in \Gamma$ let $f_b(t)$ be the unique raked trigonometric polynomial of degree $2k-1$ with 
constant term 1 that has a root at $t=a$ of multiplicity $m_a$ and a root at $t=b$ 
of multiplicity $m_b$. Let $x, y, z \in \Gamma$ be distinct points such that 
$y$ lies between $a$ and $z$ and $x$ lies between $a$ and $y$.
Suppose that $f_{y}(t) \geq 0$ for all $t \in {\Bbb S}$ and that $f_{y}$ has a root 
(of necessarily even multiplicity) in the arc $\Gamma +\pi$. 
Then $f_{x}(t)$ is positive for all $t \in {\Bbb S} \setminus \{a, x\}$
while $f_{z}(t)$ is negative for some $t \in {\Bbb S}$.
\endroster
\endproclaim 

Let us denote for a moment the maximum possible length of an arc $\Gamma$ satisfying 
Part (1) of Theorem 5.1 by $\psi_k$. In Lemma 5.7 below we prove that $\psi_k=\phi_k$, 
the maximum length of an arc with the neighborliness property of  Theorem 1.1.

\example{(5.2) Example} 

Suppose that $k=2$. The only possible set of multiplicities in 
Part (2) of Theorem 5.1 is $m_a=2$ and $m_b=2$. 
The polynomial $f(t)=1-\cos 3t$ has roots at 
$t=\pm 2\pi/3$ and a root at $t=0$, all of multiplicity $2$, while remaining non-negative on ${\Bbb S}$. 
Combining Parts (3) and (2) of Theorem 5.1 we conclude that  
$$\psi_2= {2 \pi \over 3} \approx 2.094395103.$$

Suppose that $k=3$. The only possible set of multiplicities in Part (2) of Theorem~5.1 
is $m_a =2$ and $m_b=4$. The polynomial $f(t)=1-\cos 5t$ 
has roots at $t=0, \pm 2\pi/5$, and $t=\pm 4\pi/5$, all 
of multiplicity $2$, while remaining non-negative on ${\Bbb S}$. 
Applying to $f(t)$ the deformation of Lemma 3.6 with 
$\lambda = 1/ \cos(\pi/5)$ results in the polynomial 
$f_{\lambda}(t)$ that has a root of multiplicity $4$ at $t=\pi$,
roots of multiplicity 2 at the points $\pm \alpha$ such that 
$$\cos \alpha = {\cos(2\pi/5) \over \cos (\pi/5)}={3 -\sqrt{5} \over 2},$$
and no other roots. Hence $f_{\lambda}(t)$ does not change its sign on ${\Bbb S}$.
Scaling $f_{\lambda}$, if necessary, to make the constant term 1, 
we ensure that $f_{\lambda}(t)$ is non-negative on ${\Bbb S}$. It follows by Theorem 5.1
that 
$$\psi_3 = \pi - \alpha =\pi - \arccos {3- \sqrt{5} \over 2} \approx 1.962719003.$$

Suppose that $k=4$. 
There are two possibilities for multiplicities $m_a$ and $m_b$ in Part (2) of Theorem 5.1. 
We have either $m_a=2$ and $m_b=6$ or $m_a=m_b=4$. 
It turns out that the arc satisfying the latter conditions is shorter. 
As follows from Proposition 7.6 below, we have $\psi_4 =2 \alpha$, 
where $\alpha >0$ is the smallest positive root of the equation 
$$\cos \alpha +1 -{1 \over 2} \tan^2 \alpha + 
{3 \over 8} \tan^4 \alpha - {5 \over 16} \tan^6 \alpha=0.$$
Computations show that 
$$\split \psi_4 =&2 \arccos\left( -{1 \over 48} \left(91 + 336 \sqrt{15}\right)^{1/3} + {119 \over 48\left(91 +336 \sqrt{15}\right)^{1/3}}
+{29 \over 48}\right)\\
\approx &1.870658532.\endsplit$$
In this case, the raked trigonometric polynomial $f$ of degree $7$ that has roots 
of multiplicity 4 at $t=\pm \psi_4/2$, 
also has a root of multiplicity $2$ at $t=\pi$.
\endexample 
In general, our computations suggest that in Part (2) of Theorem 5.1 
one should always choose $m_a=m_b=k$ if $k$
is even and $m_a=k+1$ and $m_b=k-1$ if $k$ is odd, but we have been unable to prove that.

To prove Theorem 5.1, we need some technical results on convergence of trigonometric polynomials.

\subhead (5.3) Convergence of trigonometric polynomials \endsubhead 
All raked trigonometric polynomials  (2.1.1) of degree at most $2k-1$ form a 
real $(2k+1)$-dimensional vector space, which we make into a normed space by letting 
$$\|f\| = \max_{t \in {\Bbb S}} |f(t)|$$
for a trigonometric polynomial $f$. For a complex polynomial $p$ of degree at most 
$4k-2$ we define 
$$\|p\| =\max_{z: \ |z|=1} |p(z)| = \max_{z:\ \|z| \leq 1} |p(z)|,$$
where the last equality follows by the maximum modulus principle for holomorphic functions.
We note that 
$$\|\PP(f)\| =\|f\|$$
for any trigonometric polynomial $f$.
We define the convergence of trigonometric and complex polynomials 
with respect to the norm $\| \cdot \|$. 

\proclaim{(5.4) Lemma} Fix a positive integer $m$. 
For a positive integer $j$, let $A_j \subset {\Bbb S}$ be a non-empty 
closed set and let $f_j(t)$ be a trigonometric polynomial of degree at most
$2k-1$ that has at least $m$ roots, counting 
multiplicities, in $A_j$. Suppose that $A_{j+1} \subset A_j$ for all 
$j$, and let 
$$B= \bigcap_{j=1}^{\infty} A_j.$$
Suppose further that for some trigonometric polynomial $f$ we have 
$$f=\lim_{j \longrightarrow +\infty} f_j.$$
Then $f$ has at least $m$ roots, counting multiplicities, in $B$. 

Suppose, in addition, that $f \not\equiv 0$, $m=2k$, $B$ 
lies in an open semicircle, and that for every $j$ the multiplicities of all roots of $f_j$
in $A_j$ are even. Then the multiplicities of all roots of $f$ in $B$ are even. 
\endproclaim 
\demo{Proof} 
Let $p_j=\PP(f_j)$. By Lemma 3.2, $p_j$ is a complex polynomial that can be written as 
$$p_j(z) = \left(z-z_{1j}\right) \cdots \left(z-z_{mj}\right) q_j(z), \tag5.4.1$$
where $q_j(z)$ is a complex polynomial of degree at most $4k-2-m$ and 
$z_{1j}, \ldots, z_{mj}$ are not necessarily distinct complex numbers 
of modulus 1 whose arguments lie in $A_j$. In addition,
$$\lim_{j \longrightarrow +\infty} p_j =p,$$
where $p(z) =\PP(f)$. From this and (5.4.1) we infer that the numbers 
$$\max_{z:\ |z|={1 \over 2}} \left| q_j(z)\right|$$
are uniformly bounded from above. Since all norms on the finite-dimensional space of 
complex polynomials of degree at most $4k-2$ are equivalent and since 
$\max_{z:\ |z|={1 \over 2}} \left| q_j(z)\right|$ is also a norm, 
it follows that the norms $\|q_j\|$ are uniformly bounded from above. (We consider the circle 
$|z|=1/2$ instead of $|z|=1$ to make sure that the factors 
$z-z_{1j}, \ldots z-z_{mj}$ in (5.4.1) are all separated from $0$.)
Hence we can find a subsequence $\left\{j_n\right\}$ such that 
$$\lim_{n \longrightarrow +\infty} q_{j_n} =q$$
for some complex polynomial $q$ and 
$$\lim_{n \longrightarrow +\infty} z_{ij_{n}} =
z^{\ast}_i, \quad \text{where} \quad z_i^{\ast} \in B \quad 
\text{for} \quad i=1, \ldots, m.$$
Then, necessarily 
$$p(z) =\left(z-z^{\ast}_1\right) \ldots \left(z-z^{\ast}_m\right) q(z).$$
Hence by Lemma 3.2 the raked trigonometric polynomial $f(t)$ has at least
 $m$ roots in $B$, counting multiplicities.
If $m=2k$ and $p \not\equiv 0$, Part (3) of Theorem~3.1 implies that 
$z^{\ast}_1, \ldots, z^{\ast}_m$ are the only roots of $p(z)$ in $B$. 
The result follows.
{\hfill \hfill \hfill} \qed
\enddemo

The following lemma plays the crucial role in our proof of Theorem 5.1.

\proclaim{(5.5) Lemma} Let $\Gamma \subset {\Bbb S}$ be an 
open arc with the endpoints $a$ and $b$ and let 
$\overline{\Gamma}$ be its closure. Let
$t_2, \ldots, t_n \in {\Bbb S} \setminus \overline{\Gamma}$ be distinct points such that the set 
$\overline{\Gamma} \cup \left\{t_2, \ldots, t_n \right\}$
lies in an open semicircle, and let $m_1, \ldots, m_n$ be positive even integers such that 
$$\sum_{i=1}^n m_i=2k.$$
For $s \in \overline{\Gamma}$, let $f_s(t)$ be the unique raked 
trigonometric polynomial of degree $2k-1$ with constant 
term 1 that has a root of multiplicity $m_i$ at $t_i$ for $i=2, \ldots, n$ 
and a root of multiplicity $m_1$ at $t=s$. If both $f_a(t)$ and $f_b(t)$ 
are non-negative on ${\Bbb S}$, then for every $s \in \Gamma$, 
the trigonometric polynomial $f_s(t)$ is positive on 
${\Bbb S} \setminus \left\{s, t_2, \ldots, t_n \right\}$.
\endproclaim
\demo{Proof} Let us consider 
$$g_s(t) ={\partial \over \partial s} f_s(t)$$ 
as in Theorem 4.2. By Lemma 4.5, 
for all $s \in \Gamma$, the point $t_i$ is a root of $g_s(t)$ of multiplicity 
at least $m_i$ for $i=2, \ldots, n$ and $s$ is a root of $g_s(t)$ of multiplicity at least $m_1-1$. 
Let ${\Bbb S}_+$ be an open semicircle containing $\overline{\Gamma}$ and the points $t_2, \ldots, t_n$.

Seeking a contradiction, let us assume that $f_{t_1}\left(t^{\ast}\right) =  0$ for some 
$t_1 \in \Gamma$ and some $t^{\ast} \in {\Bbb S} \setminus \left\{t_1, t_2, \ldots, t_n \right\}$. 
By Part (4) of Theorem 3.1, $t^{\ast} \in {\Bbb S}_+ + \pi$.
We have $f_a\left(t^{\ast}\right)  \geq 0$ and $f_b\left(t^{\ast}\right) \geq 0$. 
Therefore, the function 
$$s \longmapsto f_s\left(t^{\ast}\right)$$ attains a local minimum 
in $\Gamma$ at some point $s^{\ast}$.
Then 
$$g_{s^{\ast}}\left(t^{\ast}\right)=0 \quad \text{and} \quad f_{s^{\ast}}(t^{\ast}) \leq 0.$$
Since $f_{s^{\ast}}(t)$ has a constant term of 1, we obtain 
$$f_{s^{\ast}}(t) +f_{s^{\ast}}(t+\pi)=2 \quad \text{for all} \quad t \in {\Bbb S},$$ 
and hence 
$$t^{\ast} +\pi  \ne s^{\ast},  t_2, \ldots, t_n.$$
Since the constant term of $g_{s^{\ast}}(t)$ is 0, Part (2) of Theorem 3.1 implies that 
$g_{s^{\ast}}(t) \equiv 0$. This however contradicts Theorem 4.2.

Hence for every $s \in \Gamma$ the trigonometric polynomial $f_s(t)$ has no roots other than 
$s, t_2, \ldots, t_n$. 
By Remark 2.3(2), 
we have $f_s(t) >0$ for all $t \in {\Bbb S} \setminus\left\{s, t_2, \ldots, t_n \right\}$.
{\hfill \hfill \hfill} \qed
\enddemo

\subhead (5.6) Proof of Theorem 5.1 \endsubhead
To prove Part (1), let us choose a point $t^{\ast} \in {\Bbb S}$ and 
let us assume, seeking a contradiction, that there is a nested sequence 
of open arcs
$$\Gamma_1 \supset \Gamma_2 \supset \ldots \supset \Gamma_i \supset \ldots \tag5.6.1$$
such that 
$$\bigcap_{j=1}^{\infty} \Gamma_j = \left\{t^{\ast}\right\},$$
and such that for every $j$ there is a raked trigonometric polynomial 
$f_j(t)$ of degree $2k-1$, with 
constant term 1, with $2k$ roots, counting multiplicities, in $\Gamma_j$ and a root 
somewhere else on the circle. By Part (4) of Theorem 3.1, that additional root must lie in 
$\Gamma_j + \pi$. 
Let $h_j(t)$ be the scaling of $f_j$ to a trigonometric polynomial of norm 1. 
Then there is a subsequence of the sequence
$h_j(t)$ converging to a raked trigonometric polynomial $h$. In particular, $\|h\| =1$, and 
hence $h(t) \not\equiv 0$. 
It follows from Lemma 5.4  that $t^{\ast}$ is a root of $h$ of multiplicity at least 
$2k$ and that $t^{\ast} +\pi$ is a root of $h$. Since both $t^{\ast}$ and $t^{\ast}+\pi$ 
are roots of $h(t)$, we obtain that $h(t)$ has a zero 
constant term and that $t^{\ast} +\pi$ is, in fact, a 
root of $h(t)$ of multiplicity at least $2k$. Hence 
Part~(1) of Theorem 3.1 implies that $h(t) \equiv 0$, 
which is a contradiction.

By Remark 2.3(2), a trigonometric polynomial with constant term 1
that does not change its sign on ${\Bbb S}$ is non-negative on ${\Bbb S}$.
Finally, the example of polynomial $1-\cos (2k-1) t$ shows that the length of an
arc $\Gamma$ in Part (1) is less than $\pi$.

To prove Part (2), 
we construct a nested sequence of open arcs (5.6.1) such that 
$$\bigcap_{j=1}^{\infty} \Gamma_j = \overline{\Gamma},$$
where $\overline{\Gamma}$ is the closure of $\Gamma$. By our assumption, 
for every $j$ there is a raked trigonometric 
polynomial $f_j(t)$ of degree at most $2k-1$ that has $2k$ roots counting multiplicity 
in $\Gamma_j$ and a root elsewhere, necessarily 
in $\Gamma_j + \pi$. As in the proof of Part~(1), let us scale $f_j(t)$ to a
 trigonometric polynomial $h_j(t)$ such that $\|h_j\|=1$ 
and construct the limit trigonometric polynomial $h$. 
Then $h \not\equiv 0$, and by Lemma~5.4, $h$ has 
roots $t_1, \ldots, t_n \in \overline{\Gamma}$ of even multiplicities $m_1, \ldots, m_n$ such that 
$m_1 + \ldots + m_n =2k$, and a root $t^{\ast} \in \overline{\Gamma} + \pi$. 
By Part (2) of Theorem 3.1, $h$ has a non-zero constant term.

We rescale $h$ to a raked trigonometric polynomial $f(t)$ with constant 
term 1. Then each $t_i$ is a root of $f(t)$ of multiplicity $m_i$ and $f(t^{\ast})=0$.

Our assumption that $\Gamma$ is of maximum possible length implies 
that the endpoints $a$ and $b$ of $\Gamma$ are roots of $f(t)$.
Our goal is to show that for every $i=1, \ldots, n$ we have 
either $t_i=a$ or $t_i=b$, that is, that there are no 
roots inside $\Gamma$.

Seeking a contradiction, let us assume that $t_1 \in \Gamma$. 
We choose a closed arc $A \subset \Gamma$ with the endpoints $x$ and $y$, containing 
$t_1$ in its interior and such that $t_i \notin A$ for $i =2, \ldots, n$. For $s \in A$, let $f_s(t)$ 
be the raked trigonometric polynomial of Theorem 4.2 
that has a root at $t=s$ of multiplicity $m_1$ and a root 
at $t_i$ of multiplicity $m_i$ for $i=2, \ldots, n$. In particular, 
$$f_s =f \quad \text{if} \quad s=t_1.$$
We observe that 
$$f_s(t) \geq 0 \quad \text{for all} \quad t \in {\Bbb S} 
\quad \text{and all} \quad s \in A.$$
Indeed, if $f_s(t_0) <0$ for some $t_0 \in {\Bbb S}$ then a trigonometric polynomial 
$\hat{f}$ with constant term 1 that has a root of multiplicity $m_1$ at $s$ and roots of multiplicity
 $m_i$ at some points $\hat{t}_i \in \Gamma$ sufficiently close to $t_i$
will also satisfy $\hat{f}(t_0)<0$, which contradicts the definition of $\Gamma$. 
Hence $f_x(t) \geq 0$ for all $t \in {\Bbb S}$ and $f_y(t) \geq 0$ 
for all $t \in {\Bbb S}$. Lemma 5.5 then implies 
that $f\left(t^{\ast}\right)=f_{t_1}\left(t^{\ast}\right)>0$, which is contradiction.

To Prove Part (3), we note that for any $b \in \Gamma$ 
sufficiently close to $a$, by Part (1) of the theorem we have 
$f_{b}(t) > 0$ for all $t \in {\Bbb S} \setminus \{a, b\}$. 
We can choose such a point $b$ so that $x$ lies between
$b$ and $y$ and then $f_{x}(t) > 0$ for all 
$t \in {\Bbb S} \setminus \{a, x \}$ by Lemma 5.5. Assume now that 
$f_{z}(t) \geq 0$ for all $t \in {\Bbb S}$. Then by Lemma 5.5 we have $f_{y}(t) > 0$ for all 
$t \in {\Bbb S} \setminus \{a, y\}$, which is a contradiction.
{\hfill \hfill \hfill}\qed

\proclaim{(5.7) Lemma} Let $\psi_k$ be the maximum length of an open arc $\Gamma$ 
in Theorem 5.1 and let $\phi_k$ be the maximum length of an open arc $\Gamma$ in Theorem 1.1.
Then $\psi_k=\phi_k$.
\endproclaim
\demo{Proof} From Remark 2.3(1) it follows immediately that $\phi_k \geq \psi_k$. 

Let $\Gamma \subset {\Bbb S}$ be an open arc of length $\psi_k$ with the endpoints $a$ and $b$ and let 
$\tilde{\Gamma} \supset \Gamma$ be a closed arc  with the endpoints 
$a$ and $c$ strictly containing $\Gamma$ and lying in an open semicircle.  By Part (2) 
of Theorem 5.1 there exist positive even integers $m_a$ and $m_b$ such that $m_a+m_b=2k$ and a 
raked trigonometric polynomial $f(t)$ of degree $2k-1$ and with constant term $1$ that has a root at
$t=a$ of multiplicity $m_a$, a root at $t=b$ of multiplicity $m_b$, and some other root 
$t^{\ast} \in \Gamma +\pi$.
For $s \in \tilde{\Gamma}$ let $f_s(t)$ be the unique raked trigonometric polynomial 
of degree $2k-1$ with constant term 1 that has 
a root of multiplicity $m_a$ at $t=a$ and a root of multiplicity $m_b$ at $t=s$. 
Seeking a contradiction, let us assume 
that for any distinct $t_1, \ldots, t_k \in \tilde{\Gamma}$, the unique 
raked trigonometric polynomial of degree $2k-1$ and with constant term 1 that has roots of 
multiplicity two at $t_1, \ldots, t_k$ remains non-negative on the entire circle ${\Bbb S}$. 
As in Section 5.6, using the limit argument, we conclude that $f_c(t) \geq 0$ for all  $t \in {\Bbb S}$. 
This, however, contradicts Part (3) of Theorem 5.1 since $f_b$ has a root in $\Gamma +\pi$.

In view of Remark 2.3(1), it follows that for some distinct
 $t_1, \ldots, t_k \in \tilde{\Gamma}$, the convex hull 
$$\conv\Bigl( U\left(t_1\right), \ldots, U\left(t_k\right)\Bigr)$$
is not a face of $\BB_k$. Hence $\phi_k \leq \psi_k$.
{\hfill \hfill \hfill} \qed
\enddemo

\head 6. Neighborliness of the symmetric moment curve \endhead

In this section we prove Theorems 1.1 and 1.3. Our proofs are based on the 
following main result. 

\proclaim{(6.1) Theorem} 
For every positive integer $k$ there exists a number $\pi> \phi_k > {\pi/2}$ such that if 
$\Gamma \subset {\Bbb S}$ is an open arc of length $\phi_k$, if $t_1, \ldots, t_n \in \Gamma$ 
are distinct points, and $m_1, \ldots, m_n$ are positive even integers such that 
$$\sum_{i=1}^n m_i =2k,$$
then the unique raked trigonometric polynomial $f(t)$ of degree $2k-1$ with constant term 1 that 
has a root of multiplicity $m_i$ at $t_i$ for $i=1, \ldots, n$ is 
positive everywhere else on the circle ${\Bbb S}$.
\endproclaim

The proof is based on Theorem 5.1 and the following lemma.

\proclaim{(6.2) Lemma} 
Let $f(t)$ be the raked trigonometric polynomial of degree $2k-1$ with constant term 
$1$ that has a root of multiplicity $2m$ at $t=0$ 
and a root of multiplicity $2n$ at $t=\pi/2$, where $m$ and $n$ are 
positive integers such that $m+n =k$. Then $f(t)$ has no other roots in the circle ${\Bbb S}$.
\endproclaim 
\demo{Proof} 
We have 
$$f(t)=1+\sum_{j=1}^k a_j \cos (2j-1) t + \sum_{j=1}^k b_j \sin (2j-1) t$$
for some real $a_j$ and $b_j$.
In addition,
$$f'(0) = \ldots = f^{(2m-1)}(0) = 0 \quad \text{and} \quad 
f'(\pi/2) = \ldots = f^{(2n-1)}(\pi/2)=0. \tag6.2.1$$

Let 
$$a(t)=\sum_{j=1}^k a_j \cos (2j-1) t \quad \text{and} \quad b(t) = \sum_{j=1}^k b_j \sin (2j-1) t,$$
so that
$$f(t) =1+a(t) +b(t) \quad \text{and} \quad f'(t) =a'(t) +b'(t). \tag6.2.2$$
Observe that
$${d^{2r-1} \over dt^{2r-1}} a(t) \Big|_{t=0}=0 \quad \text{and} 
\quad  {d^{2r} \over dt^{2r}} b(t) \Big|_{t=0} =0 \tag6.2.3$$
for any positive integer $r$, and 
$${d^{2r} \over dt^{2r}} a(t) \Big|_{t=\pi/2} =0 \quad \text{and} \quad 
{d^{2r-1} \over dt^{2r-1}} b(t) \Big|_{t=\pi/2} =0 \tag6.2.4$$
for any positive integer $r$.

Combining (6.2.1) -- (6.2.4) we conclude that $t=0$ is a root of $a'(t)$ 
of multiplicity at least $2m-1$ and a root of 
$b'(t)$ of multiplicity at least $2m$. Similarly, $t=\pi/2$ 
is a root of $a'(t)$ of multiplicity at least $2n$ and a root of 
$b'(t)$ of multiplicity at least $2n-1$. Since $f(0)=0$, we obtain that 
$a(t) \not\equiv 0$, and hence $a'(t)\not\equiv 0$.
Also since $f(\pi/2)=0$, it follows that $b(t) \not\equiv 0$, 
and hence $b'(t) \not\equiv 0$. By Part (2) of Theorem~3.1,
the trigonometric polynomial $a'(t)$ has a root of multiplicity 
$2m-1$ at $t=0$, a root of multiplicity $2n$ at $t=\pi/2$ 
and no other roots in the circle ${\Bbb S}$, 
while the trigonometric polynomial $b'(t)$ has a root of multiplicity $2m$ at $t=0$, 
a root of multiplicity $2n-1$ at $t=\pi/2$ and no other roots in the circle.

We conclude that the functions $a(t)$ and $b(t)$ are monotone on the interval $0 < t < \pi/2$. 
Since $a(0)=-1$ and $a(\pi/2)=0$, we infer that $a(t)$ is monotone increasing for $0 < t < \pi/2$, 
and hence $a(t) < 0$ for all $0 < t < \pi/2$. Since $b(0)=0$ and $b(\pi/2)=-1$, 
we obtain that $b(t)$ is monotone 
decreasing for $0 < t < \pi/2$, and therefore $b(t) < 0$ for all $0 < t < \pi/2$. 
As
$$a(t+\pi) = -a(t) \quad \text{and} \quad b(t+\pi)=-b(t),$$
it follows that $a(t) >0$ for $\pi < t < 3\pi/2$ and $b(t) >0$ for $\pi < t < 3\pi/2$. Therefore, 
$$f(t) \geq 1 \quad \text{for all} \quad \pi \leq t \leq 3\pi/2. $$
The latter equation yields the result, as by Part (4) of Theorem 3.1, 
a root $t^{\ast}$ of $f(t)$ distinct from $0$ and $\pi/2$, if exists, must 
satisfy $\pi \leq t^{\ast} \leq 3\pi/2$.
{\hfill \hfill \hfill} \qed
\enddemo

\subhead (6.3) Proof of Theorem 6.1 \endsubhead 
 By Part (1) of Theorem 5.1, there exists a number $\eta_k >0$ such that if a raked trigonometric 
polynomial $f(t)$ of degree $2k-1$ and with a constant term 1 has roots at $t=0$ and $t=\eta_k$ with 
positive even multiplicities summing up to $2k$, then $f(t)$ is positive everywhere else. 
It follows from Lemma 6.2 and Lemma 5.5 that the same remains true for all
 $0 < \eta_k \leq \pi/2$. Using the shift $f(t) \longmapsto f(t+a)$ of 
raked trigonometric polynomials, we conclude that for every arc 
$\Gamma \subset {\Bbb S}$ of length not exceeding 
$\pi/2$, a raked trigonometric polynomial $f(t)$ of degree $2k-1$ 
with constant term 1 that has roots of even 
multiplicities summing up to $2k$ at the endpoints of $\Gamma$ 
remains positive everywhere else in ${\Bbb S}$.
The proof now follows from Part (2) of Theorem 5.1.
{\hfill \hfill \hfill} \qed

\subhead (6.4) Proofs of Theorem 1.1 and 1.3 \endsubhead 
Theorem 1.1 follows from Theorem 6.1 and Remark 2.3(1), while Theorem 1.3 follows from 
Remark 2.3(1) and Lemma 5.5.
{\hfill \hfill \hfill} \qed 

\head 7. The limit of neighborliness \endhead

In this section, we prove Theorem 1.2. 
Our goal is to construct a raked trigonometric polynomial $f_k(t)$ of degree 
$2k-1$ such that $f_k(t)$ has a root of multiplicity $2k-2$ at $t=0$, 
roots of multiplicity $2$ each at $t=\pm \beta_k$ 
for some $\pi/2 < \beta_k < \pi$, and such that 
$f_k(t) \geq 0$ for all $t \in {\Bbb S}$. 
It then follows from Theorem 5.1 and Lemma 5.7
that $\phi_k \leq \beta_k$, and establishing that 
$\beta_k \longrightarrow \pi/2$ as $k$ grows, we complete the proof.

\proclaim{(7.1) Lemma} The function 
$$f(t)=\sin^{2k-1}t$$ is a raked trigonometric polynomial of degree $2k-1$.
\endproclaim 
\demo{Proof} We have 
$$\split \sin^{2k-1} t = &\left({e^{it} - e^{-it} \over 2i} \right)^{2k-1}=
{1 \over (-4)^{k-1}} {1 \over 2i} \sum_{j=0}^{2k-1} {2k-1 \choose j} (-1)^{j} e^{i(2k-2j-1)t} \\
=&{1 \over (-4)^{k-1}}\sum_{j=0}^{k-1} {2k-1 \choose j}  
 \left( {(-1)^j e^{i(2k-2j-1)t} + (-1)^{2k-1-j}e^{i(2j-2k+1)t} \over 2i}\right) \\
 =&{1 \over (-4)^{k-1}} \sum_{j=0}^{k-1} {2k-1 \choose j} (-1)^j  \sin (2k-2j-1) t.  \endsplit$$
{\hfill \hfill \hfill} \qed
\enddemo

\proclaim{(7.2) Lemma} For $k\geq 1$ let 
$$h_k(t)=\int_0^t \sin^{2k-1}(\tau) \ d \tau.$$ 
Then $h_k(t)$ is a raked trigonometric polynomial of 
degree $2k-1$ and $t=0$ is a root of $h_k(t)$ of multiplicity $2k$. Moreover,
$$h_k(t)={(2k-2)!! \over (2k-1)!!} 
\left(1-(\cos t) \sum_{j=0}^{k-1} {(2j-1)!! \over (2j)!!} \sin^{2j} t \right),$$
where we agree that $0!!=(-1)!!=1$.
\endproclaim
\demo{Proof} From Lemma 7.1, 
$h_k(t)$ is a raked trigonometric polynomial of degree $2k-1$. Moreover, $h_k(0)=0$ and 
$h_k'(t)=\sin^{2k-1}t$, from which it follows that $t=0$ is a root of $h_k(t)$ of multiplicity $2k$. 
Since 
$$h_1(t)=\int_0^t \sin \tau \ d\tau = 1-\cos t,$$
and since
for $n >1$,  
$$\int_0^t \sin^n \tau \ d\tau =
 -{1 \over n} \left(\sin^{n-1} t\right)\left( \cos t\right) 
+ {n-1 \over n} \int_0^t \sin^{n-2} \tau \ d\tau,$$
we obtain by induction that
$$\int_0^t \sin^{2k-1} \tau \ d\tau =
{(2k-2)!! \over (2k-1)!!}\left(1 - \cos t \sum_{j=0}^{k-1} {(2j-1)!! \over (2j)!!} \sin^{2j} t\right),$$
as claimed.
{\hfill \hfill \hfill} \qed
\enddemo

\proclaim{(7.3) Lemma} Let $h_k(t)$ be the trigonometric polynomial defined in Lemma 7.2 and let 
$$F_k(t) = \sin^2(t) h_{k-1}(t)-h_k(t).$$
Then there exists a unique 
$${\pi \over 2} \ < \ \beta_k \ < \ \pi$$
 such that 
$$F_k\left(\beta_k\right)=0.$$
In addition, 
$$\lim_{k \longrightarrow +\infty} \beta_k = {\pi \over 2}.$$
\endproclaim 
\demo{Proof} From Lemma 7.2, we deduce
$$\split &F_k\left({\pi \over 2}\right)= h_{k-1}\left({\pi \over 2}\right)-h_k\left({\pi \over 2}\right)
={(2k-4)!! \over (2k-3)!!} -{(2k-2)!! \over (2k-1)!!} \ > \  0 \quad \text{and} \\ 
&F_k(\pi)=-h_k(\pi)=-2{(2k-2)!! \over (2k-1)!!} \ < \ 0. \endsplit$$
Moreover,
$$F_k'(t)=2(\sin t )(\cos t) h_{k-1}(t)-h_k'(t)+ \sin^2 (t) h_{k-1}'(t)=2(\sin t)(\cos t)h_{k-1}(t).$$
In particular, $F_k'(t)<0$ for $\pi/2 < t < \pi$, and hence  
$F_k(t)$ is decreasing on the interval $\pi/2 < t < \pi$. 
Since $F_k(\pi/2)>0$ and $F_k(\pi)<0$,
there is a unique $\pi/2 < \beta_k < \pi$ such that $F_k\left(\beta_k\right)=0$.

To find the limit behavior of $\beta_k$, we use the expansion
$$(1-x)^{-1/2} =\sum_{j=0}^{\infty} {(2j-1)!! \over (2j)!!} x^j \quad \text{for real} \quad -1 < x < 1.$$
Substituting $x=\sin ^2 t$ we obtain 
$$\sum_{j=0}^{\infty} {(2j-1)!! \over (2j)!!} \sin^{2j} t =-{1 \over  \cos t} \quad \text{provided} \quad \pi/2 < t < \pi.$$
Hence from Lemma 7.2, for $\pi/2 < t < \pi$ we have 
$$\split h_k(t)=&{(2k-2)!! \over (2k-1)!!} \left(1 -(\cos t) \left( \sum_{j=0}^{\infty} {(2j-1)!! \over (2j)!!} \sin^{2j} t -
\sum_{j=k}^{\infty} {(2j-1)!! \over (2j)!!} \sin^{2j} t  \right)\right)\\=&{(2k-2)!! \over (2k-1)!!}\left(2 +
(\cos t) \sum_{j=k}^{\infty}{(2j-1)!! \over (2j)!!} \sin^{2j} t \right)  \endsplit$$
and
$$ \split &{(2k-3)!! \over (2k-4)!!} F_k(t)=\\ &\quad 2 \sin^2 t -2 {2k-2 \over 2k-1} + (\cos t) \sum_{j=k}^{\infty} 
\left( {(2j-3)!! \over (2j-2)!!}-{2k-2 \over 2k-1} 
{(2j-1)!! \over (2j)!!} \right) \sin^{2j} t.\endsplit$$
It follows that $F_k(t) <0$ for every $\pi/2 < t < \pi$ such that $\sin^2 t \leq (2k-2)/(2k-1)$. Since $F_k(t)$ is decreasing
for $\pi/2 < t < \pi$, we conclude that 
$$\sin^2 \beta_k \ > \ {2k-2 \over 2k-1}, \tag7.3.1$$ 
and hence 
$$\lim_{k \longrightarrow +\infty} \beta_k ={\pi \over 2},$$
as desired.
{\hfill \hfill \hfill} \qed
\enddemo

\proclaim{(7.4) Lemma} 
Let $h_k(t)$ be the trigonometric polynomial defined in Lemma 7.2 and let $\beta_k$ 
be the number defined in Lemma 7.3. 
Let 
$$f_k(t) =\sin^2 \left( \beta_k \right) h_{k-1}(t)-h_k(t).$$
Then $f_k(t)$ is a raked trigonometric polynomial of degree $2k-1$ 
such that $t=0$ is a root of $f_k(t)$ of multiplicity 
$2k-2$, $t=\pm \beta_k$ are the roots of multiplicity 2 each 
and $f_k(t) \geq 0$ for all $t \in {\Bbb S}$.
\endproclaim
\demo{Proof} It follows by Lemma 7.2 that $f_k(t)$ is a 
trigonometric polynomial of degree $2k-1$ and that 
$t=0$ is a root of $f_k(t)$ of multiplicity at least $2k-2$. 
From the definition of $\beta_k$ in Lemma 7.3, 
we conclude that $t=\beta_k$ is a root of $f_k(t)$. Moreover, since
$$f_k'(t)=\sin^{2k-1} t -\left(\sin^2 \beta_k\right) \sin^{2k-3} t,$$
we have $f'\left(\beta_k\right)=0$, so
the multiplicity of the root at $t=\beta_k$ is at least 2. 
By Part (3) of Theorem 3.1, the multiplicities of the roots 
at $t=0$ and $t=\beta_k$ are $2k-2$ and 2 respectively and 
there are no other roots of $f_k(t)$ in the open arc 
$0 < t < \pi$. Also, by Lemma 7.2 and (7.3.1),  we have
$$f_k(\pi)=2 \sin^2 \left(\beta_k \right) {(2k-4)!! \over (2k-3)!!} - 2 {(2k-2)!! \over (2k-1)!!} >0.$$
Since $f_k(-t)=f_k(t)$, we conclude that $t=-\beta_k$ is a root of 
$f_k(t) \geq 0$ of multiplicity $2$ and that 
$f_k(t)>0$ for all $t \ne 0, \pm \beta_k$.
{\hfill \hfill \hfill} \qed
\enddemo

\subhead (7.5) 
Proof of Theorem 1.2 \endsubhead 
Let $\psi_k$ be the maximum length of an open arc satisfying 
Part (1) of Theorem 5.1. It follows from Lemma 7.4 that $\psi_k \leq \beta_k$, 
and hence from Lemma 5.7 that $\phi_k \leq \beta_k$.
Lemma 7.3 then yields the proof.
{\hfill \hfill \hfill} \qed
\bigskip

All available computational evidence suggests that for even $k$ 
the smallest length of the arc in Part (2) of Theorem 5.1 
is achieved when the multiplicites $m_a$ and $m_b$ are equal: $m_a=m_b=k$. 
The following results provides an 
explicit equation for the length of such an arc.

\proclaim{(7.6) Proposition} 
Suppose that $k$ is even. Let $\alpha_k>0$ be the smallest number such that 
the necessarily unique raked trigonometric polynomial $f(t)$ of degree $2k-1$ 
with constant term 1 that has roots at 
$t=\alpha_k$ and $t=-\alpha_k$ of multiplicity $k$ each also has a root $t^{\ast}$ 
elsewhere in ${\Bbb S}$. 
Then $t^{\ast}=\pi$ and $\alpha_k$ is the smallest positive root of the equation $F(\alpha)=0$ where
$$F(\alpha) = \cos \alpha + 1 + \sum_{j=1}^{k-1} (-1)^j {(2j-1)!! \over (2j)!!} 
\tan^{2j} \alpha. \tag7.6.1$$
\endproclaim
\demo{Proof}
We note that the raked trigonometric polynomial $\tilde{f}(t)=f(-t)$ also has a root of multiplicity 
$k$ at $t=\alpha_k$ and a root of multiplicity $k$ at $t=-\alpha_k$. By Corollary 3.5, we must have 
$\tilde{f}(t)=f(t)$, and hence 
$$f(t)=1 + \sum_{j=1}^k a_j \cos (2j-1)t$$
for some real $a_1, \ldots, a_k$. Then the raked trigonometric polynomial $f'(t)$ has  roots 
at $t=\alpha_k$, $-\alpha_k$, $\alpha_k +\pi$, and $-\alpha_k +\pi$ of multiplicity $k-1$ each 
as well as roots at $t=0$ and $t=\pi$.
By Part (1) of Theorem 3.1, $f'(t)$ has no other roots and $a_k\ne 0$. 
By Part (4) of Theorem 3.1, the root $t^{\ast}$ must lie
in an open arc $\Gamma$ with the endpoints $\alpha_k +\pi$ and $-\alpha_k+\pi$. 
From the definition of $\alpha_k$, it follows
that $f(t) \geq 0$ for all $t \in \Gamma$, and hence $t^{\ast}$ is a local minimum of $f(t)$. Thus 
$f'\left(t^{\ast}\right)=0$, and so $t^{\ast}=\pi$. 
Moreover, $t^{\ast} =\pi$ is a root of $f(t)$ of multiplicity 2.

We choose 
$$\lambda = {1 \over \cos \alpha_k}$$
in Lemma 3.6 and consider the $\lambda$-deformation $f_{\lambda}(t)$ of $f(t)$. 
Let 
$$p=\PP(f) \quad  \text{and}\quad  p_{\lambda}=\PP\left(f_{\lambda}\right).$$
Since $\alpha_k$ and $-\alpha_k$ are roots of $f(t)$ of multiplicity $k$ each, 
the complex numbers $e^{i \alpha_k}$ and 
$e^{-i \alpha_k}$ are roots of $p$ of multiplicity $k$ each. 
Then $z =1$ is a root of $p_{\lambda}(z)$ of multiplicity 
$2k$, and hence $t=0$ is a root of $f_{\lambda}$ of multiplicity $2k$. 

As $t=\pi$ is a root of $f$ of multiplicity $2$, it follows that 
$z=-1$ is a root of $p(z)$ of multiplicity $2$. 
Thus, 
$${-1 +\sin \alpha_k \over \cos \alpha_k} \quad \text{and} \quad 
{-1 - \sin \alpha_k \over \cos \alpha_k} \tag7.6.2$$
are roots of $p_{\lambda}(z)$.

Since $t=0$ is a root of multiplicity $2k$ of $f_{\lambda}(t)$ the trigonometric polynomial $f_{\lambda}(t)$ should be 
proportional 
to the trigonometric polynomial $h_k(t)$ of Lemma 7.2. Therefore,
$$\aligned p_{\lambda}(z) = &\gamma z^{2k-1}\left(1 - \left({z +z^{-1} \over 2}\right) 
\left(1 + \sum_{j=1}^{k-1}  {(2j-1)!! \over (2j)!!} \left( {z -z^{-1} \over 2i} \right)^{2j}
\right) \right) \\ &\text{for some} \quad \gamma \ne 0. \endaligned \tag7.6.3$$
Substituting either of the roots of (7.6.2) in (7.6.3), we obtain the desired equation 
$$F\left(\alpha_k\right)=0.$$

Suppose now that some number $0 < \alpha < \pi/2$ also satisfies the equation $F(\alpha)=0$. Then 
$${-1 +\sin \alpha \over \cos \alpha} \quad \text{and} \quad  
{-1-\sin \alpha \over \cos \alpha} \tag7.6.4$$
are roots of polynomial $q=\PP\left(h_k\right)$, 
where $h_k(t)$ is the trigonometric polynomial of Lemma 7.2.
Let us choose $\lambda = \cos \alpha$ and let $g_{\lambda}(t)$ 
be the $\lambda$-deformation of $h_k(t)$ as 
in Lemma 3.6. Let $q_{\lambda}=\PP\left(g_{\lambda}\right)$. 
Since the numbers introduced in (7.6.4) are roots of $q$, we conclude 
that $z=-1$ is a root of multiplicity 2 of $q_{\lambda}(z)$, 
and hence $t=\pi$ is a root of multiplicity $2$ of $g_{\lambda}(t)$.
Similarly, since $t=0$ is a root of multiplicity $2k$ of $h_k(t)$, 
we conclude that $z=1$ is a root of multiplicity $2k$ 
of $q(z)$, and hence the numbers $e^{i \alpha}$ and $e^{-i \alpha}$ 
are roots of  $q_{\lambda}$, each of multiplicity $k$. 
Therefore, $t=\alpha$ and $t=-\alpha$ are roots of $g_{\lambda}(t)$, 
each of  multiplicity $k$. It then follows,
by minimality of $\alpha_k$, that $\alpha \geq \alpha_k$, which completes the proof.
{\hfill \hfill \hfill} \qed
\enddemo

\head Acknowledgment \endhead

The authors are grateful to anonymous referees for careful reading of the paper and helpful suggestions.

\enddocument